\documentclass[11pt,oneside,american]{amsart}
\usepackage{mathptmx}
\usepackage[T1]{fontenc}
\usepackage[latin9]{inputenc}
\usepackage{geometry}
\usepackage{fancyhdr}
\usepackage{amsthm}

\makeatletter
\numberwithin{equation}{section}
\numberwithin{figure}{section}
\theoremstyle{plain}
\newtheorem{thm}{\protect\theoremname}
\theoremstyle{plain}
\newtheorem{lem}[thm]{\protect\lemmaname}

\makeatother

\usepackage{babel}
\providecommand{\lemmaname}{Lemma}
\providecommand{\theoremname}{Theorem}

\begin{document}
\title{fast convolution algorithm for state space models.}
\author{Gregory Beylkin}
\address{Department of Applied Mathematics, University of Colorado at Boulder,
UCB 526, Boulder, CO 80309-0526}
\begin{abstract}
We present an unconditionally stable algorithm for applying matrix
transfer function of a linear time invariant system (LTI) in time
domain. The state matrix of an LTI system used for modeling long range
dependencies in state space models (SSMs) has eigenvalues close to
$1$. The standard recursion defining LTI system becomes unstable
if the $m\times m$ state matrix has just one eigenvalue with absolute
value even slightly greater than 1. This may occur when approximating
a state matrix by a structured matrix to reduce the cost of matrix-vector
multiplication from $\mathcal{O}\left(m^{2}\right)$ to $\mathcal{O}\left(m\right)$
or $\mathcal{O}\left(m\log m\right).$ We introduce an unconditionally
stable algorithm that uses an approximation of the rational transfer
function in the z-domain by a matrix polynomial of degree $2^{N+1}-1$,
where $N$ is chosen to achieve any user-selected accuracy. Using
a cascade implementation in time domain, applying such transfer function
to compute $L$ states requires no more than $2L$ matrix-vector multiplications
(whereas the standard recursion requires $L$ matrix-vector multiplications).
However, using unconditionally stable algorithm, it is not necessary
to assure that an approximate state matrix has all eigenvalues with
absolute values strictly less than 1 i.e., within the desired accuracy,
the absolute value of some eigenvalues may possibly exceed $1$. Consequently,
this algorithm allows one to use a wider variety of structured approximations
to reduce the cost of matrix-vector multiplication and we briefly
describe several of them to be used for this purpose.
\end{abstract}

\maketitle

\section{Introduction}

In \cite{GU-GO-RE:2021,GU:2023,GU-DAO:2024} authors describe an approach
to sequence modeling using state space models (SSMs) to address long
range dependencies. The setup of SSM uses a linear multiple-input
multiple-output (MIMO) time invariant system (LTI), 
\begin{eqnarray}
x'\left(t\right) & = & Ax\left(t\right)+Bu\left(t\right)\nonumber \\
y\left(t\right) & = & Cx\left(t\right)+Du\left(t\right),\label{eq:LTI}
\end{eqnarray}
where $A$ is an $m\times m$ state matrix, $B$, $C$ and $D$ are
$m\times p$, $q\times m$ and $q\times p$ matrices with $1\le p,q\le m$.
Here $x$, $u$ and $y$ are vector functions of dimension $m$, $p$
and $q$. Typically $p$ and $q$ are significantly smaller than $m$,
e.g. in \cite{GU-GO-RE:2021} $p=q=1$. A discretization of the ordinary
differential equation in (\ref{eq:LTI}) yields a system of the form
\begin{eqnarray}
\mathbf{x}_{n} & = & \overline{A}\mathbf{x}_{n-1}+\overline{B}\mathbf{u}_{n}\nonumber \\
\mathbf{y}_{n} & = & C\mathbf{x}_{n}+D\mathbf{u}_{n},\label{eq:Discrete=000020LTI}
\end{eqnarray}
where, using for example the trapezoidal rule with step size $\Delta$,
we have $\overline{A}=\left(I-\frac{\Delta}{2}A\right)^{-1}\left(I+\frac{\Delta}{2}A\right)$
and $\overline{B}=\Delta\left(I-\frac{\Delta}{2}A\right)^{-1}B$.
Alternative discretization schemes for the ordinary differential equation
in (\ref{eq:LTI}) can be used with step size $\Delta$ playing the
role of time scale in SSM schemes (see \cite{GU-GO-RE:2021,GU-DAO:2024,GU:2023}).
We assume that the magnitude of the eigenvalues of the matrix $\overline{A}$
is less than $1$.

Modeling long range dependencies leads to a dense matrix $\overline{A}$
with absolute values of eigenvalues close to $1$. As reported in
\cite{GU-GO-RE:2021}, using an example of the so-called HiPPO matrix
(see Section~\ref{sec:An-example-of=000020HiPPO}), diagonalization
of the matrix $\overline{A}$ yields ill-conditioned eigenvectors.
In \cite{GU-GO-RE:2021,GU-DAO:2024,GU:2023} the matrix $\overline{A}$
is approximated by a structured matrix, i.e. by changing basis it
becomes diagonal-plus-low-rank (see also \cite{Y-N-M-M-E:2023}).
Alternatively, using the Fourier transform as a way to accelerate
application of a convolution was suggested in \cite{GU-GO-RE:2021}
as well as in \cite{F-E-N-T-Z-D-R-R:2023}. Importantly, using the
recurrence (\ref{eq:Discrete=000020LTI}), it is necessary to assure
that structured approximation yields a matrix with magnitude of all
eigenvalues less than $1.$

Given an input sequence $\left\{ \mathbf{u}_{n}\right\} _{n=0}^{\infty}$,
the goal is to obtain the output sequence $\left\{ \mathbf{y}_{n}\right\} _{n=0}^{\infty}$,
which can be accomplished by convolution with a transfer function.
Using the z-transform (described below for completeness), the transfer
function between these two sequences is obtained as
\begin{equation}
H\left(z\right)=C\left(I_{m}-z^{-1}\overline{A}\right)^{-1}\overline{B}+D,\label{eq:transfer=000020function}
\end{equation}
where $I_{m}$ is the $m\times m$ identity matrix. The question is
how to apply efficiently the operator $H$ as a convolution in time
domain and avoid issues associated with possible accumulation of error
when using the recurrence (\ref{eq:Discrete=000020LTI}) directly.

In this paper we propose a simple unconditionally stable algorithm
of applying the operator $H$ in time domain. Additional improvements
using structured matrices can be used within the new algorithm to
further reduce the computational cost, as briefly described in Section~\ref{sec:Structured-representation-of}.

\section{Fast cascade algorithm}

In \cite{BEYLKI:1995,BE-LE-MO:2012} we introduced a novel method
of accurate approximation of causal and anti-causal IIR filters by
a cascade of FIR filters leading to automatic design of efficient
digital filters. Many properties, such as filters with exact linear
phase or symmetric filters satisfying the perfect reconstruction condition,
can only be obtained using non-causal IIR filters. Our method of approximating
causal and anti-causal IIR filters produces FIR filters that, within
any user-selected accuracy, inherit properties of IIR filters. The
resulting cascades of FIR filters have a straightforward parallel
implementation. A simple technical tool underpinning these results
leads to a fast algorithm for applying the operator $H$ in (\ref{eq:transfer=000020function})
in time domain. For completeness, we present a detailed derivation.

Using the z-transforms 
\[
X\left(z\right)=\sum_{n=0}^{\infty}\mathbf{x}_{n}z^{-n},\,\,\,\,U\left(z\right)=\sum_{n=0}^{\infty}\mathbf{u}_{n}z^{-n}\,\,\,\,\mbox{and}\,\,\,Y\left(z\right)=\sum_{n=0}^{\infty}\mathbf{y}_{n}z^{-n},
\]
we obtain from (\ref{eq:Discrete=000020LTI}) (with $\mathbf{x}_{-1}=\mathbf{0}$)
\[
\begin{cases}
X\left(z\right)=z^{-1}\overline{A}X\left(z\right)+\overline{B}U\left(z\right)\\
Y\left(z\right)=CX\left(z\right)+DU\left(z\right)
\end{cases}
\]
yielding $X\left(z\right)=\left(I_{m}-z^{-1}\overline{A}\right)^{-1}\overline{B}U\left(z\right)$
and
\[
Y\left(z\right)=H\left(z\right)U\left(z\right),
\]
where $H$ is given in (\ref{eq:transfer=000020function}).

Since the magnitude of all eigenvalues of $\overline{A}$ is less
than $1$, Lemma~\ref{lem:Let--be} below implies
\[
\left(I_{m}-z^{-1}\overline{A}\right)^{-1}=\prod_{n=0}^{\infty}\left[I_{m}+\left(z^{-1}\overline{A}\right)^{2^{n}}\right],
\]
where $\left|z\right|=1$ and 
\begin{equation}
H\left(z\right)=C\prod_{n=0}^{\infty}\left[I_{m}+\left(z^{-1}\overline{A}\right)^{2^{n}}\right]\overline{B}+D.\label{eq:cascade=000020polynomial}
\end{equation}
Approximating $H$ by a product with $N+1$ terms , we define
\begin{equation}
H_{N}\left(z\right)=C\prod_{n=0}^{N}\left[I_{m}+\left(z^{-1}\overline{A}\right)^{2^{n}}\right]\overline{B}+D.\label{eq:cascade=000020approximation}
\end{equation}
Effectively, we are approximating a rational matrix function $\left(I_{m}-z^{-1}\overline{A}\right)^{-1}$
by a matrix polynomial. If the magnitude of some eigenvalues of the
matrix $A$ is close to $1$, then the degree $2^{N+1}-1$ of the
matrix polynomial $H_{N}$ can be large. However, the cost of applying
this matrix polynomial as a cascade in time domain requires only $N+1$
matrix-vector multiplications. We have
\begin{lem}
\label{lem:Let--be}Let $\gamma<1$ be largest singular value of $\overline{A}$
and $N$ a positive integer. Since $\gamma\ge\left|\lambda_{j}\right|$,
$j=1,\dots,m$, we estimate
\[
\left\Vert \left(I_{m}-z^{-1}\overline{A}\right)^{-1}-\prod_{n=0}^{N}\left(I_{m}+\left(z^{-1}\overline{A}\right)^{2^{n}}\right)\right\Vert \leq\frac{\left\Vert \overline{A}\right\Vert ^{2^{N+1}}}{1-\left\Vert \overline{A}\right\Vert }=\frac{\gamma^{2^{N+1}}}{1-\gamma}
\]
and, as a result, we have an approximation for $\left|z\right|=1$,
\[
\left\Vert H\left(z\right)-H_{N}\left(z\right)\right\Vert \le\frac{\gamma^{2^{N+1}}}{1-\gamma}\left\Vert C\right\Vert \left\Vert \overline{B}\right\Vert ,
\]
where
\begin{equation}
H_{N}\left(z\right)=C\prod_{n=0}^{N}\left(I_{m}+\left(z^{-1}\overline{A}\right)^{2^{n}}\right)\overline{B}+D\label{eq:convolution=000020operator}
\end{equation}
By taking $N$ large enough, we can achieve any user-selected accuracy.
\end{lem}

Proof of Lemma~\ref{lem:Let--be} uses induction on $N$ and is similar
to the proof of convergence for the Neumann series, see \cite{BEYLKI:1995,BE-LE-MO:2012}. 

In order to apply $\left(I_{m}-z^{-1}\overline{A}\right)^{-1}$ in
time domain, we need to compute the powers of the matrix $\overline{A}$,
$\overline{A}^{2},\overline{A}^{4},\overline{A}^{8},\dots$ .

\section{Convolution via a cascade algorithm\protect\label{sec:Applying-convolution-via}}

Given an input sequence $\left\{ \mathbf{u}_{\ell}\right\} _{\ell=0}^{L-1},$
the representation of the matrix polynomial in (\ref{eq:convolution=000020operator})
yields a fast cascade algorithm for applying convolution to obtain
the sequence $\left\{ \mathbf{y}_{\ell}\right\} _{\ell=0}^{L-1}$.
Computing the sequence $\mathbf{y}_{\ell}$ using (\ref{eq:Discrete=000020LTI})
and setting $\mathbf{x}_{-1}=0$ yields
\begin{equation}
\mathbf{y}_{\ell}=\sum_{k=0}^{\ell}C\overline{A}^{\ell-k}\overline{B}\mathbf{u}_{k}+D\mathbf{u}_{k},\label{eq:sum}
\end{equation}
an expression that uses all powers of matrix $\overline{A}$. Lemma~\ref{lem:Let--be}
allows us to reuse intermediate matrix-vector products by the matrix
$\overline{A}$ and, as a result, we need to apply only the powers
$\overline{A}^{2^{n}}$, $n=0,\dots,N-1$. Note that the powers of
$z^{-2^{n}}$ in (\ref{eq:convolution=000020operator}) correspond
to a shift by $2^{n}$ in time domain. 

We start algorithm by computing 
\[
\mathbf{v}_{\ell}^{\left(0\right)}=\overline{B}\mathbf{u}_{\ell},\,\,\,\,\ell=0,1,\dots,L-1
\]
and forming an $m\times L$ matrix
\[
V^{\left(0\right)}=\left[\mathbf{v}_{0}^{\left(0\right)}\mathbf{v}_{1}^{\left(0\right)}\dots\mathbf{v}_{L-1}^{\left(0\right)}\right].
\]
Next we update columns of $V^{\left(0\right)}$ starting from column
$\ell=1$ according to
\begin{equation}
\begin{cases}
\mathbf{v}_{\ell}^{\left(1\right)}=\overline{A}\mathbf{v}_{\ell-1}^{\left(0\right)}+\mathbf{v}_{\ell}^{\left(0\right)} & 1\le\ell\le L-1\\
\mathbf{v}_{\ell}^{\left(1\right)}=\mathbf{v}_{\ell}^{\left(0\right)} & \ell=0
\end{cases}\label{eq:step1}
\end{equation}
and, for convenience of notation, rename the matrix as $V^{\left(1\right)}$.
We then update columns of $V^{\left(1\right)}$starting from $\ell=2$,
\begin{equation}
\begin{cases}
\mathbf{v}_{\ell}^{\left(2\right)}=\overline{A}^{2}\mathbf{v}_{\ell-2}^{\left(1\right)}+\mathbf{v}_{\ell}^{\left(1\right)} & 2\le\ell\le L-1\\
\mathbf{\mathbf{v}_{\ell}^{\left(2\right)}=v}_{\ell}^{\left(1\right)} & 0\le\ell\le1
\end{cases}\label{eq:step2}
\end{equation}
and rename the result as $V^{\left(2\right)}$. Next we update columns
of $V^{\left(2\right)}$starting from $\ell=4$,
\begin{equation}
\begin{cases}
\mathbf{v}_{\ell}^{\left(3\right)}=\overline{A}^{4}\mathbf{v}_{\ell-4}^{\left(2\right)}+\mathbf{v}_{\ell}^{\left(2\right)} & 4\le\ell\le L-1\\
\mathbf{\mathbf{v}_{\ell}^{\left(3\right)}=v}_{\ell}^{\left(2\right)} & 0\le\ell\le3,
\end{cases}\label{eq:step3}
\end{equation}
and so on. At each step $\ell=0,2,\dots L-1$ we update columns of
$V^{\left(n-1\right)}$ in the range $2^{n-1}\le\ell\le L-1$ to yield
$V^{\left(n\right)}$,
\[
\begin{cases}
\mathbf{v}_{\ell}^{\left(n\right)}=\overline{A}^{2^{n-1}}\mathbf{v}_{\ell-2^{n-1}}^{\left(n-1\right)}+\mathbf{v}_{\ell}^{\left(n-1\right)} & 2^{n-1}\le\ell\le L-1\\
\mathbf{v}_{\ell}^{\left(n\right)}=\mathbf{v}_{\ell}^{\left(n-1\right)} & 0\le\ell\le2^{n-1}-1.
\end{cases}
\]
Finally, we compute\textbf{
\[
\mathbf{y}_{\ell}=C\mathbf{v}_{\ell}^{\left(N\right)}+D\mathbf{u}_{\ell},\,\,\,\ell=0,1,\dots L-1.
\]
}As described, the algorithm requires $\mathcal{O}\left(mL\right)$
memory and $\mathcal{O}\left(m^{2}\log L+mL\right)$ arithmetic operations.
Using approximation of matrix $\overline{A}$ by a structured matrix,
the cost of arithmetic operations can be further reduced to $\mathcal{O}\left(m\log L+mL\right)$
or $\mathcal{O}\left(m\log m\log L+mL\right)$ depending on the selected
representation of matrix $\overline{A}$, see Section~\ref{sec:Structured-representation-of}.
Noting that sizes of matrices $\overline{B}$, $C$ and $D$, $m\times p$,
$q\times m$ and $q\times p$, are typically much smaller than the
$m\times m$ size of $\overline{A}$, matrix-vector multiplications
by matrices $\overline{B}$, $C$ and $D$ are significantly less
expensive and we did not count them in the complexity estimates. Our
complexity estimate is for the case $p=q=1$, a single-input single-output
(SISO) model as in \cite{GU-GO-RE:2021}.

While the degree of matrix polynomial $H_{N}$ in (\ref{eq:convolution=000020operator})
is $2^{N+1}-1$ , the cascade algorithm requires only $N+1$ matrix-vector
multiplications. For example, if $N=14$ then the degree of the matrix
polynomial is $2^{N+1}-1=32767$, whereas applying it as a convolution
requires only 15 matrix-vector multiplications. While providing an
approximation with any finite user-selected accuracy, the cascade
algorithm is unconditionally stable in contrast with that using recurrence
in (\ref{eq:Discrete=000020LTI}) directly.

\section{\protect\label{sec:An-example-of=000020HiPPO}An example of HiPPO
matrix}

The following triangular $m\times m$ matrix is used in \cite{GU-GO-RE:2021,G-J-T-R-R:2022,GU:2023,Y-N-M-M-E:2023},

\begin{equation}
A=-\begin{cases}
\sqrt{2n+1}\sqrt{2k+1} & n>k\\
n+1 & n=k\\
0 & n<k,
\end{cases}\label{eq:hippo=000020matrix}
\end{equation}
where $1\le n,k\le m$. This is a so-called HiPPO (High-order Polynomial
Projection Operator) matrix, see \cite{G-D-E-R-R:2020}. Using it
as an example and setting $m=100$, we compute 
\begin{equation}
\overline{A}=\left(I-\frac{\Delta}{2}A\right)^{-1}\left(I+\frac{\Delta}{2}A\right),\label{eq:hippo-discrete}
\end{equation}
with $\Delta=0.5\cdot10^{-3}$. The diagonal elements of the resulting
triangular matrix, i.e. its eigenvalues, range from $d_{1}=0.999000499750125$
to $d_{100}=0.9507437210436478$. Since $d_{1}^{2^{15}}\approx5.87548*10^{-15}$,
the matrix $\overline{A}^{2^{15}}$ can be neglected and, as a result,
$15$ matrix-vector multiplications allow us to apply the matrix polynomial
$H_{14}$ of degree $32767$. Using double precision arithmetic in
Mathematica$^{TM}$, the powers of $\overline{A}$ are computed with
accuracy $\approx10^{-14}$. We note that the powers of the matrix
$\overline{A}$ can always be computed with a large number of accurate
digits and used with as many digits as an application requires.

\section{\protect\label{sec:Structured-representation-of}Structured approximations
of matrices}

If a structured approximation of $m\times m$ matrix $\overline{A}$
is available, then the cost of matrix-vector multiplications using
powers of this matrix can be significantly reduced. For example, by
changing basis, the matrix $\overline{A}$ is modified to be diagonal-plus-low-rank
in \cite{GU-GO-RE:2021}.

We observe that there are other structured approximations that can
be used to change complexity of matrix-vector multiplication from
$\mathcal{O}\left(m^{2}\right)$ to $\mathcal{O}\left(m\right)$ or
$\mathcal{O}\left(m\log m\right)$.

For example, the matrix in (\ref{eq:hippo-discrete}) admits the so-called
partitioned low rank representation (PLR) (see e.g. \cite{BEY-SAN:2005}),
where the rank of the off-diagonal blocks is equal to $1$. The cost
of applying the $m\times m$ matrix (\ref{eq:hippo-discrete}) and
its powers in PLR representation does not exceed $\mathcal{O}\left(m\log m\right)$.

Other structured approximations include a wavelet representation of
the matrices $A$ or $\overline{A}$, where an approximation via the
so-called non-standard form \cite{BE-CO-RO:1991} yields an $\mathcal{O}\left(m\right)$
algorithm for matrix-vector multiplication. It is possible to use
multiwavelets for the same purpose, see \cite{ALPERT:1993}. 

The choice of a structured representation also depends on the selected
discretization of the LTI system in (\ref{eq:LTI}). For example,
for the HiPPO matrix, if one selects $\overline{A}=\exp\left(\Delta A\right)$
and $\overline{B}=\left(\Delta A\right)^{-1}\left[\exp\left(\Delta A\right)-I\right]\Delta A$,
then a wavelet representation of $\exp\left(\Delta A\right)$ is more
efficient than that for $\overline{A}=\left(I-\frac{\Delta}{2}A\right)^{-1}\left(I+\frac{\Delta}{2}A\right)$. 

A class of matrices that admit a sparse structured approximation with
complexity of matrix-vector product $\mathcal{O}\left(m\right)$ or
$\mathcal{O}\left(m\log m\right)$ is large and such matrices can
potentially be considered in setting up the LTI system. If a MIMO
system is used, then structured approximations can also be used for
matrices $B$ or $\overline{B}$ and $C$ and $D$.

\section{Discussion}

The algorithm described in this paper is formally slower than the
direct use of the recurrence (\ref{eq:Discrete=000020LTI}) by a factor
not exceeding two ($\overline{A}$ is applied to every entry of the
sequence, $\overline{A}^{2}$ to every second entry, $\overline{A}^{4}$
to every forth entry, etc.). This raises a question ``what are the
advantages of replacing the recurrence (\ref{eq:Discrete=000020LTI})
with a cascade convolution in Section~\ref{sec:Applying-convolution-via}?''.
If the absolute values of all eigenvalues of the matrix $\overline{A}$
in (\ref{eq:Discrete=000020LTI}) are significantly smaller than $1$,
then it is likely that the eigenvalues of a structured matrix constructed
as its approximation will also be smaller than $1$. In such cases
the recurrence (\ref{eq:Discrete=000020LTI}) is stable and it is
not necessary to use the algorithm in Section~\ref{sec:Applying-convolution-via}.
However, a structured approximation of a matrix $\overline{A}$ with
absolute values of eigenvalues close to $1$ may end up with eigenvalues
that violate this condition (within the accuracy of approximation).
In such cases the recurrence in (\ref{eq:Discrete=000020LTI}) is
unstable whereas the cascade convolution in Section~\ref{sec:Applying-convolution-via}
is always stable. Note that in applications with a goal of analyzing
long range dependencies, at least some eigenvalues of a matrix $\overline{A}$
will have absolute value close to $1$. In such problems, with a modest
increase in computational cost, the algorithm in Section~\ref{sec:Applying-convolution-via}
provides for a wider choice of possible structured approximations
and allows one to consider alternative bases than those used in \cite{GU-GO-RE:2021,G-J-T-R-R:2022,GU:2023,Y-N-M-M-E:2023}.

\section{Acknowledgment}

Part of this research was performed while the author was visiting
the Institute for Pure and Applied Mathematics (IPAM), which is supported
by the National Science Foundation (Grant No. DMS-1925919).

\bibliographystyle{ieeetr}

\end{document}